\documentclass[12pt]{article}
\usepackage{amsfonts,amssymb,epsfig}
\usepackage{amsmath}

\date{}
\begin{document}
\newtheorem{df}{Definition}
\newtheorem{thm}{Theorem}
\newtheorem{lm}{Lemma}
\newtheorem{pr}{Proposition}
\newtheorem{co}{Corollary}
\newtheorem{re}{Remark}
\newtheorem{note}{Note}
\newtheorem{claim}{Claim}
\newtheorem{problem}{Problem}
\newtheorem{fact}{Fact}

\def\R{{\mathbb R}}

\def\E{\mathbb{E}}
\def\calF{{\cal F}}
\def\N{\mathbb{N}}
\def\calN{{\cal N}}
\def\calH{{\cal H}}
\def\n{\nu}
\def\a{\alpha}
\def\d{\delta}
\def\t{\theta}
\def\e{\varepsilon}
\def\t{\theta}
\def\pf{ \noindent {\bf Proof: \  }}
\def\trace{\rm trace}
\newcommand{\qed}{\hfill\vrule height6pt
width6pt depth0pt}
\def\endpf{\qed \medskip} \def\colon{{:}\;}
\setcounter{footnote}{0}

\def\Lip{{\rm Lip}}

\renewcommand{\qed}{\hfill\vrule height6pt  width6pt depth0pt}

\title{Three observations regarding Schatten $p$ classes
}

\author{Gideon Schechtman\thanks{Supported in part by the Israel Science Foundation. } } \maketitle

\begin{abstract}
The paper contains three results, the common feature of which is that they deal with the Schatten $p$ class. The first is a presentation of a new complemented subspace of $C_p$ in the reflexive range (and $p\not= 2$). This construction answers a question of Arazy and Lindestrauss from 1975. The second result relates to tight embeddings of finite dimensional subspaces of $C_p$ in $C_p^n$ with small $n$ and shows that $\ell_p^k$ nicely embeds into $C_p^n$ only if $n$ is at least proportional to $k$ (and then of course the dimension of $C_p^n$ is at least of order $k^2$). The third result concerns single elements of $C_p^n$ and shows that for $p>2$ any $n\times n$ matrix of $C_p$ norm one and zero diagonal admits, for every $\e>0$, a $k$-paving of $C_p$ norm at most $\e$ with $k$ depending on $\e$ and $p$ only.
\end{abstract}

\section{Introduction}
\subsection{Complemented subspaces}
Recall that for $1\le p<\infty$ $C_p$ denotes the Banach space of all compact operators $A$ on $\ell_2$ for which the norm $\|A\|_p=(\trace(A^*A)^{p/2})^{1/p}$ is finite. Determining the complemented subspaces of these spaces was a subject of investigation for quite a while. In particular Arazy and Lindenstrauss in \cite{al} list nine isomorphically distinct infinite dimensional complemented subspaces of $C_p$, $1<p\not= 2<\infty$. They are all complemented by quite natural projections which are given by setting certain subsets of the entries of a given matrix to zero. They are also all norm one projections. Unlike the situation with $L_p$ there is some hope of isomorphically characterizing all complemented subspaces of $C_p$. The main result of section \ref{sec:complemented} of this paper is to introduce a new complemented subspace of $C_p$, solving some problems from \cite{al}. To describe the space we need some more notations.
For a finite or infinite matrix $A$ we denote by $A(k,l)$ its $k,l$ entry. $Z_p$ will denote the Banach space of matrices whose norm
\[
\|A\|_{Z_p}=(\sum_{k=1}^\infty(\sum_{l=1}^\infty|A(k,l)|^2)^{p/2})^{1/p}
\]
is finite.
For $p>2$,  $\widetilde{Z_p}$ will denote the Banach space of matrices whose norm
\[
\|A\|_{\widetilde{Z_p}}=(\|A\|_{Z_p}^p+ \|A^*\|_{Z_p}^p)^{1/p}
\]
is finite. For $1\le q<2$, $\widetilde{Z_q}$ will denote the Banach space of matrices whose norm
\[
\|A\|_{\widetilde{Z_q}}=\inf\{(\|B\|_{Z_q}^q+ \|C^*\|_{Z_q}^q)^{1/q}\ ; \ A=B+C\}
\]
is finite. Note that if $q=p/(p-1)$, $p>2$, then $\widetilde{Z_q}$ is the dual of $\widetilde{Z_p}$ and vice versa.\\
We also denote by $C_p^n, Z_p^n$ and $\widetilde{Z_p^n}$ the spaces of $n\times n$ matrices with the norms inherited from $C_p, Z_p$ and $\widetilde{Z_p}$ respectively.
Define
\[
D_p=(\bigoplus_{n=1}^\infty\widetilde{Z_p^n})_p.
\]

The main result in section \ref{sec:complemented} is that $D_p$ is complemented in $C_p$, $1<p<\infty$, and is not isomorphic to any of the previously known complemented subspaces of $C_p$.\\
The proof depends heavily on a result from \cite{lp} showing that $\widetilde{Z_p}$ is the ``unconditional version" of $C_p$.

\subsection{Tight Embeddings} Before describing the results of section \ref{sec:tight} we would like to motivate it by describing what is known in the $L_p$ case. \\
Given a $k$-dimensional subspace $X$ of $L_p(0,1)$ one can ask what is the minimal $n$ such that $X$ embeds with constant $2$, say, into $\ell_p^n$. This was extensively studied and, up to $\log$ factors basically solved. Except for a power of $\log k$, $n$ can be taken to be of order $k$ for $1\le p<2$ and of order $k^{p/2}$ for $2<p<\infty$. These orders are best possible, again up to the $\log$ factor, as is seen by looking at $X=\ell_2^k$. See \cite{js} for a survey of these results. In the case that $p$ is an even integer the $\log$ factor is not needed (\cite{s}). \\
It is natural to seek similar results in $C_p$. Given a $k$-dimensional subspace $X$ of $C_p$ what is the smallest $n$ such that $X$ $2$-embeds into $C_p^n$? If $X$ is $\ell_2^k$ the answer is known for quite a while: $n\approx \sqrt k$ for $1\le p<2$ and $n\approx_p k^{\frac{p}{p+2}}$ for $2<p<\infty$ \cite{flm}. One could guess that also here the worst case is $\ell_2^k$. In section \ref{sec:tight} we show that this is wrong and $\ell_p^k$ is worse than $\ell_2^k$ in this respect.

We prove in Theorem \ref{thm:tight} that if $C_p^n$ contains a $2$-isomorphic copy of $\ell_p^k$ then $n\gtrsim_p k$. This holds for all $1\le p\not= 2<\infty$.\\
The proof is algebraic in nature, estimating the rank of some operator.

\subsection{Paving}
Recall that a $k$-paving of a $n\times n$ matrix $A$ is a matrix of the form $\sum_{i=1}^k P_{\sigma_i}AP_{\sigma_i}$ where $\{\sigma_1,\dots,\sigma_k\}$ is a partition of $\{1,\dots,n\}$ and for a subset $\sigma\subseteq \{1,\dots,n\}$ $P_{\sigma}AP_{\sigma}$ is the matrix whose $k,l$ term is $A(k,l)$ if both $k$ and $l$ are in $\sigma$ and $0$ otherwise.

Marcus, Spielman and Srivastava \cite{mss} solved recently the Kadison--Singer conjecture which is equivalent to the paving conjecture. In our terms they proved:

\medskip

\noindent{\em For each $\e>0$ there is a $k=k(\e)$ such that for all $n$ and all $n\times n$ matrix $A$ with zero diagonal and norm one (as an operator on $\ell_2^n$) there is a $k$-paving of norm at most $\e$.}

\medskip

Given a reasonable norm on the space of $n\times n$ matrices, one can ask if a similar result holds for that norm. In \cite{bhkw} it was proved that this is the case for self adjoint matrices and the norms $C_4$ and $C_6$ (the case of $C_2$ is easy and was known before). The hope of the authors of \cite{bhkw} was that one will be able to prove a similar result for all zero diagonal matrices for a sequence of $C_{p_n}$ norms with $p_n\to\infty$ and where the function $k(\e)$ is independent of $n$. Then the paving conjecture would follow.

In Theorem \ref{thm:paving} of section \ref{sec:paving} we show that the paving conjecture indeed holds for all the $C_p$ norms, $2<p<\infty$. However, the function $k$ that we get depends on $p$.\\
The proof in \cite{bhkw} (for $p=4,6$) is by an averaging argument. Our proof is also quite a simple averaging argument. The difference with the argument in \cite{bhkw} is that instead of estimating the average of $\|P_\sigma A P_\sigma\|$ over all possible $\sigma$ from above we estimate the average of $\|A- P_\sigma A P_\sigma\|$ from below and then use the uniform convexity of $C_p$.

\section{A new complemented subspace of $C_p$, $1<p\not=2<\infty$}\label{sec:complemented}

We start with stating two facts relating the spaces $C_p$, $Z_p$ and $\widetilde{Z_p}$.
The first is relatively simple:

\begin{fact}\label{fact:1}
For $p>2$, $\|A\|_p\ge\|A\|_{Z_p}$ and consequently also $2^{1/p}\|A\|_p\ge\|A\|_{\widetilde{Z_p}}$.
\end{fact}
Let us reproduce the (known) proof of this fact: Let $E_{k,l}$ denote the matrix with 1 in the $k,l$ place and zeros elsewhere. Also for $\e=(\e_1,\dots,\e_k,\dots)$, with $\e_k=\pm1$ for all $k$, denote $\Delta_\e={\rm diag}(\e_1,\dots,\e_k,\dots)$. By the unitary invariance of the $C_p$ norm and the convexity, over positive definite matrices, of $B\to\trace (B^{p/2})$, we get that for all $A=\sum_{k,l}a_{k,l}E_{k,l}$,
\begin{multline*}
\|A\|_p^p={\rm Ave}_{\e}\|\Delta_\e A\|_p^p={\rm Ave}_{\e}\trace(\Delta_\e AA^*\Delta_\e)^{p/2}\ge\trace[{\rm Ave}_{\e}\Delta_\e AA^*\Delta_\e]^{p/2}.
\end{multline*}
Now, as is easily seen,
\begin{multline*}
{\rm Ave}_{\e}\Delta_\e AA^*\Delta_\e=\sum_{k,l}(AA^*)(k,l){\rm Ave}_{\e}\e_k\e_lE_{k,l}\\=\sum_k(AA^*)(k,k)E_{k,k}=
\sum_k\sum_l|a_{k,l}|^2E_{k,k}
\end{multline*}
and
\[
\trace(\sum_k\sum_l|a_{k,l}|^2E_{k,k})^{p/2}=\|A\|_{Z_p}^p.
\]
\endpf

\bigskip
\noindent As was noted by the referee, one can also prove Fact \ref{fact:1} by interpolating between the simpler cases $p=2$ and $p=\infty$. A dual statement (extending also to $q=1$) holds for $1\le q<2$. We'll not use this here.

By $\alpha\gtrsim_p\beta$ (resp. $\alpha\lesssim_p\beta$) we mean that for some positive constant $K$, depending only on $p$, $K\alpha\ge \beta$ (resp. $\alpha\le K\beta$). $\alpha\approx_p\beta$ means  $\alpha\gtrsim_p\beta$ and $\alpha\lesssim_p\beta$. If the subscript $p$ is omitted the constant $K$ is meant to be an absolute constant.

The second Fact is much more profound and follows from a result of Lust-Piquard \cite{lp} for $1<p<\infty$ and of Lust-Piquard and Pisier \cite{lpp} for $p=1$. See also \cite{px}.

\begin{fact}\label{fact:2}
Let $\{\e_{k,l}\}_{k,l=1}^\infty$ be a matrix of independent random variables each taking the values $\pm1$ with equal probability. For a matrix $A=\sum_{k,l=1}^\infty a_{k,l}E_{k,l}$ in $C_p$, $1\le p<\infty$,
\[
({\rm Ave}_\e\|\sum_{k,l=1}^\infty \e_{k,l}a_{k,l}E_{k,l}\|_p^p)^{1/p}\approx_p\|A\|_{\widetilde{Z_p}}.
\]
\end{fact}

Let $\{R_i^n\}_{i=1}^{2^{n^2}}$ be a list of all the $n\times n$ sign matrices; i.e., all matrices each of whose entries is $1$ or $-1$.

Given two matrices of the same dimensions $A$ and $B$, we denote by $A\circ B$ their Schur product; i.e., the matrix whose $i,j$ component is the product of the $i,j$ component of $A$ with the $i,j$ component of $B$.

Given an $n\times n$ matrix $A$ let $\overline{A}$ be the $n2^{n^2}\times n2^{n^2}$ matrix given by
\[
\overline{A}(k,l)=\left\{
                    \begin{array}{ll}
                      (R_i^n\circ A)(k-(i-1)n,l-(i-1)n),\ & (i-1)n<k,l\le in, \\
                      &  i=1,\dots,2^{n^2} \\
                      0 & \mbox{otherwise} \\
                    \end{array}
                  \right.
 \]
i.e., $\overline{A}$ is a block diagonal matrix with $2^{n^2}$ blocks each of size $n\times n$ and the $i$-th block is $R_i^n\circ A$. Note that by Fact \ref{fact:2} above, for all $1\le p<\infty$,
\[
\|\overline{A}\|_p\approx_p 2^{n^2/p} \|A\|_{\widetilde{Z_p}}.
\]
Set
\[
D_p^n=\{\overline{A}\ ;\ A\in C_p^n\},\ \mbox{equipped with the norm induced by}\ C_p^{n2^{n^2}}.
\]
Then for each $1\le p<\infty$, $A\to 2^{-n^2/p}\overline{A}=\overline{2^{-n^2/p}A}$ is an isomorphism, with constant depending only on $p$, between $\widetilde{Z_p^n}$ and $D_p^n$.

\begin{pr}\label{pr:main} For $1<p<\infty$ $D_p^n$ is $\lambda_p$-complemented in $C_p^{n2^{n^2}}$ where $\lambda_p$ depends on $p$ only.
\end{pr}

\pf
For $B\in C_p^{n2^{n^2}}$ and $i=1,\dots,2^{n^2}$ let $B_i$ be the $n\times n$ central block of $B$ supported on the coordinates in $((i-1)n,in]$; i.e.,
\[
B_i(k,l)=B(k+(i-1)n,l+(i-1)n)\ \ , \ \ 1\le k,l\le n .
 \]
For $B\in C_p^{n2^{n^2}}$ let
\[
\Phi(B)=2^{-n^2}{\sum_{i=1}^{2^{n^2}}R_i^n\circ B_i}
\] and define
\[
Q(B)=2^{-n^2}\overline{\sum_{i=1}^{2^{n^2}}R_i^n\circ B_i}.
\]
Hence $Q(B)=\overline{\Phi(B)}\in D_p^n$. Since for every $A\in C_p^n$
\[
\Phi(\overline{A})=2^{-n^2}\sum_{i=1}^{2^{n^2}}R_i^n\circ R_i^n\circ A=A,
 \]
 $Q$ is obviously a projection onto $D_p^n$. It is also easy to see that it is self adjoint. It is therefore enough to evaluate its norm as an operator on $C_p$, $p>2$. As we remarked above,
\begin{equation}\label{eq:Q(B)}
\begin{array}{rl}
\|Q(B)\|_p&\approx_p 2^{n^2/p} \|2^{-n^2}\sum_{i=1}^{2^{n^2}}R_i^n\circ B_i\|_{\widetilde{Z_p}}\\
&= 2^{n^2/p}(\|2^{-n^2}\sum_{i=1}^{2^{n^2}}R_i^n\circ B_i\|_{ Z_p}^p+\|2^{-n^2}(\sum_{i=1}^{2^{n^2}}R_i^n\circ B_i)^*\|_{Z_p}^p)^{1/p}.
\end{array}
\end{equation}
We'll show that
\begin{equation}\label{eq:Z_p}
2^{n^2/p}\|2^{-n^2}\sum_{i=1}^{2^{n^2}}R_i^n\circ B_i\|_{ Z_p}\le \|B\|_p.
\end{equation}
Since the treatment of the other term on the right hand side of (\ref{eq:Q(B)}) is identical this will prove the Proposition. (\ref{eq:Z_p}) is really just the statement that $Q$ is a norm $1$ projection on $Z_p$ but let us repeat the proof.
\[\begin{array}{rl}
\|\sum_{i=1}^{2^{n^2}}R_i^n\circ B_i\|_{Z_p}^p&=\sum_{k=1}^n\|\sum_{i=1}^{2^{n^2}}(R_i^n\circ B_i)(k,\cdot)\|_{\ell_2}^p\\
&\le\sum_{k=1}^n(\sum_{i=1}^{2^{n^2}}\|(R_i^n\circ B_i)(k,\cdot)\|_{\ell_2})^p\\
&\quad \quad \mbox{by the triangle inequality,}\\
&=\sum_{k=1}^n(\sum_{i=1}^{2^{n^2}}\|B_i(k,\cdot)\|_{\ell_2})^p\\
&\le \sum_{k=1}^n2^{n^2(p-1)}\sum_{i=1}^{2^{n^2}}\|B_i(k,\cdot)\|_{\ell_2}^p\\
&\quad \quad \mbox{by Holder's inequality},\\
&\le 2^{n^2(p-1)}\sum_{k^{\prime}=1}^{{n2^{n^2}}}\|B(k^{\prime},\cdot)\|_{\ell_2}^p\\
&\quad \quad \mbox{by denoting}\ k^{\prime}=in+k,\\
&=2^{n^2(p-1)}\|B\|_{Z_p}^p\le 2^{n^2(p-1)}\|B\|^p_p\\
&\quad \quad \mbox{by Fact}\  \ref{fact:1}.
\end{array}
\]
This proves  (\ref{eq:Z_p}) and thus concludes the proof.
\endpf

Recall the notation from \cite{al}: $S_p=(\bigoplus_{n=1}^\infty C_p^n)_p$; i.e., $S_p$ is the subspace of $C_p$ spanned by disjoint diagonal blocks, the $n$-th one being of size $n\times n$. We also denote: $D_p=(\bigoplus_{n=1}^\infty \widetilde{Z_p^n})_p$. As a corollary we immediately get the main result of this section.

\begin{thm}\label{thm:main}
For $1<p\not=2<\infty$ $D_p$ is isomorphic to a complemented subspace of $S_p$ and thus of $C_p$. It is not isomorphic to any of the nine previously known infinite dimensional complemented subspaces of $C_p$ listed in Theorem 5 of \cite{al}.
\end{thm}

\pf
The first assertion follows immediately from Proposition \ref{pr:main}. Since $D_p$ has an unconditional basis and cotype $p\vee 2$ and since $C_p^n$ does not uniformly embed into a lattice with non trivial cotype with constant independent of $n$ (see \cite[Theorem 2.1] {p}), it remains to show that $D_p$ is not isomorphic to any of the spaces  listed in Theorem 5 of \cite{al} having an unconditional basis. These are the four spaces $\ell_2$, $\ell_p$, $\ell_2\oplus\ell_p$ or $Z_p$. These four spaces are isomorphic to complemented subspaces of $L_p(0,1)$, so it is enough to show that $D_p$ is not. Indeed, for $1<p<2$ $D_p$ is not even isomorphic to a subspace of $L_p(0,1)$. This follows for example from \cite{hrs}. Indeed, if the $\widetilde{Z_p^n}$-s, $1\le p<2$, uniformly embed in $L_p(0,1)$ then by a simple limiting argument so would $\widetilde{Z_p}$. However, Corollary 1.3 in \cite{hrs} states in particular that ${\rm Rad} C_p$, which is isomorphic to $\widetilde{Z_p}$, does not isomorphically embed in $L_p(0,1)$.  Since $D_p^*=D_{p/(p-1)}$ it follows that, for $p>2$, $D_p$ is not isomorphic to a complemented subspace of $L_p(0,1)$.
\endpf

\begin{rm}
1. Theorem \ref{thm:main} solves the problem in Remark (ii) on page 107 of \cite {al}. In particular, $D_p$ is complemented in $S_p$ but is not isomorphic to $S_p$ or $\ell_p$. One can of course ask whether these three spaces exhaust all isomorphic types of infinite dimensional complemented subspaces of $S_p$.\\
2. Combining the space $D_p$ with the previously known complemented subspaces, we can get two more isomorphically different complemented subspaces of $C_p$. These are $\ell_2\oplus_p D_p$ and $Z_p\oplus_p D_p$. It is not hard to show that they are not isomorphic to any of the other spaces and not isomorphic to each other. So all together we get twelve isomorphically different infinite dimensional complemented subspaces of $C_p$, $1<p\not= 2<\infty$.  This leaves open the main problem of whether there are infinitely many isomorphic classes of infinite dimensional complemented subspaces of $C_p$.\\
We also remark that for any subset $\{m_n\}_{n=1}^\infty$ of the natural numbers with $\sup m_n=\infty$, the space $(\bigoplus_{n=1}^\infty\widetilde{Z_p^{m_n}})_p$ is isomorphic to $D_p$.\\
3. Is $\widetilde{Z_p}$ isomorphic to a complemented subspace of $C_p$? We believe not. Probably, for $1\le p<2$, $\widetilde{Z_p}$ does not even isomorphically embed into $C_p$. Note that for $p>2$ the situation is different: $\widetilde{Z_p}$ isometrically embeds even in $Z_p$.
\end{rm}\\
4. It is easy to see (using \cite{lpp}) that the norm of $Q$ in the proof of Proposition \ref{pr:main} is of order $\sqrt p$ for $p>2$.

\section{Tight embeddings in $C_p^n$}\label{sec:tight}
The result of this section was obtained when I visited the University of Alberta where I greatly benefitted from interaction on the subject with Sasha Litvak and Nicole Tomczak--Jaegermann. The main result of this section is

\begin{thm}\label{thm:tight}
If $T_1,\dots,T_k$ are $n\times n$ matrices with $\|T_i\|_{C_p^n}\ge 1$ for all $i$ and
\[
{\rm Ave}_{\pm}\|\sum_{i=1}^k \pm T_i\|_{C_p^n}^p\le Kk, \ \ \mbox{if} \ \ p>2,
\]
or $\|T_i\|\le 1$ for all $i$ and
\[
{\rm Ave}_{\pm}\|\sum_{i=1}^k \pm T_i\|_{C_p^n}^p\ge K^{-1}k, \ \ \mbox {if} \ \ 1\le p<2,
\]
then $n\ge K^{\frac{-2}{|p-2|}}k$.
In particular if $T:\ell_p^k\to X\subseteq C_p^n$ is a linear isomorphism then $n\ge (\|T\|\|T^{-1}\|)^{\frac{-2p}{|p-2|}}k$.
\end{thm}

We start with two claims.

\begin{claim}\label{claim:p>2}
Let $p>2$ and let $A_i$, $i=1,2,\dots,k$, be $n\times n$ positive definite matrices with ${\trace}A_i^{p/2}\ge 1$ for all $i$. Assume ${\trace}(\sum_{i=1}^k A_i)^{p/2}\le Kk$. Then $n\ge {\rm rank}(\sum_{i=1}^k A_i)\ge K^{\frac{-2}{p-2}}k$.
\end{claim}

\pf Let $d={\rm rank}(\sum_{i=1}^k A_i)$ and let $\lambda_1,\dots,\lambda_d$ be the non zero eigenvalues of $\sum_{i=1}^k A_i$. Then

\begin{multline*}
Kk\ge {\trace}(\sum_{i=1}^k A_i)^{p/2}=\sum_{i=1}^d\lambda_i^{p/2}\ge d^{\frac{2-p}{2}}
(\sum_{i=1}^d\lambda_i)^{p/2}\ \ \mbox{by Holder's inequality}\\
=d^{\frac{2-p}{2}}(\sum_{i=1}^k {\trace}A_i)^{p/2}\ge d^{\frac{2-p}{2}}(\sum_{i=1}^k ({\trace}A_i^{p/2})^{2/p})^{p/2}\ge d^{\frac{2-p}{2}}k^{p/2}.
\end{multline*}
So  $d\ge K^{\frac{-2}{p-2}}k$.
\endpf

\begin{claim}\label{claim:p<2}
Let $0<p<2$ and let $A_i$, $i=1,2,\dots,k$, be $n\times n$ positive definite matrices with ${\trace}A_i^{p/2}\le 1$ for all $i$. Assume ${\trace}(\sum_{i=1}^k A_i)^{p/2}\ge ck$. Then $n\ge {\rm rank}(\sum_{i=1}^k A_i)\ge c^{\frac{2}{2-p}}k$.
\end{claim}

\pf Let $d={\rm rank}(\sum_{i=1}^k A_i)$ and let $\lambda_1,\dots,\lambda_d$ be the non zero eigenvalue of $\sum_{i=1}^k A_i$.

\begin{multline*}
ck\le {\trace}(\sum_{i=1}^k A_i)^{p/2}=\sum_{i=1}^d\lambda_i^{p/2}\le d^{\frac{2-p}{2}}
(\sum_{i=1}^d\lambda_i)^{p/2}\\
=d^{\frac{2-p}{2}}(\sum_{i=1}^k {\trace}A_i)^{p/2}\le d^{\frac{2-p}{2}}(\sum_{i=1}^k ({\trace}A_i^{p/2})^{2/p})^{p/2}\le d^{\frac{2-p}{2}}k^{p/2}.
\end{multline*}
So  $d\ge c^{\frac{2}{2-p}}k$.
\endpf

\noindent{\bf Proof of Theorem \ref{thm:tight}} By the easy part of the inequality in \cite{lp} (or see \cite{px}) which was actually known before and whose proof is quite easy and similar to that of Fact \ref{fact:1},
\[
{\rm Ave}_{\pm}\|\sum_{i=1}^k \pm T_i\|_{C_p^n}^p\ge {\rm tr}(\sum_{i=1}^k  T_i^*T_i)^{p/2}, \ \ \mbox{if} \ \  p>2,
\]
and
\[
{\rm Ave}_{\pm}\|\sum_{i=1}^k \pm T_i\|_{C_p^n}^p\le {\rm tr}(\sum_{i=1}^k  T_i^*T_i)^{p/2}, \ \ \mbox{if} \ \  0<p<2.
\]
Now apply Claim \ref{claim:p>2} or \ref{claim:p<2} with $A_i=T_i^*T_i$.

To prove the last claim in the statement of the theorem, assume $p>2$ and assume (as we may by multiplying T by a constant)  that $\|T^{-1}\|=1$. Letting $T_i$ be the image by $T$ of the i-th unit basis vector in $\ell_p^k$ we see that $\|T_i\|_{C_p^n}\ge 1$ and that ${\rm Ave}_{\pm}\|\sum_{i=1}^k \pm T_i\|_{C_p^n}^p\le \|T\|^p k$. Now apply the first part of the theorem with $K=\|T\|^p$. The case $1\le p<2$ is treated similarly starting with $\|T\|=1$.
\endpf

\bigskip

Is $\ell_p^k$ the worst (maybe up to $\log$ factors) $k$-dimensional subspace for tight embedding in $C_p^n$? i.e., does any $k$ dimensional subspaces of $C_p$ $2$-embed into $C_p^n$ with $n$ proportional to $k$, except maybe for a multiplicative factor of a power of $\log k$?

\section{Paving in $C_p$, $p>2$}\label{sec:paving}
The main result here is the following Theorem which clearly gives, by iteration, the result claimed in the third subsection of the Introduction.

\begin{thm}\label{thm:paving}
Let $A$ be a $2m\times 2m$ matrix with zero diagonal. Then there are mutually orthogonal diagonal
(i.e., with range a span of a subset of the natural basis) projections
$P,Q$ of rank $m$ such that,
for all $2\le p<\infty$,
\[
\|PAP+QAQ\|_p\le \left(1-\frac{1}{2^p}\right)^{1/p}\|A\|_p.
\]
\end{thm}

Given a $2m\times 2m$ matrix $A=\{a_{i,j}\}_{i,j=1}^{2m}$ and a subset $\sigma\subset\{1,2,\dots,2m\}$ of
cardinality $m$, let $A_\sigma$ denote the matrix whose $i,j$ element is 0 if $i,j$ both
belong to $\sigma$ or both belong to the complement $\sigma^c$ of $\sigma$, and $a_{i,j}$ otherwise.

\begin{pr}\label{pr:Asigma}
For all $2\le p<\infty$ and every $2m\times 2m$ matrix $A$ with zero diagonal,
\[
\rm{Ave}_{\sigma}\rm{Tr}(A_{\sigma}^*A_{\sigma})^{p/2}\ge\frac{1}{2^p}\rm{Tr}(A^*A)^{p/2},
\]
where the average is taken over all subsets of $\{1,2,\dots,2m\}$ of cardinality $m$.
\end{pr}

\pf
Applying Proposition 2 in \cite{pe} we get that
\begin{equation}\label{eq:conv}
\rm{Ave}_{\sigma}\rm{Tr}(A_{\sigma}^*A_{\sigma})^{p/2}\ge\rm{Tr}(\rm{Ave}_{\sigma} A_{\sigma}^*A_{\sigma})^{p/2}.
\end{equation}

Now, the $i,j$ element of $A_\sigma^*A_{\sigma}$ is
\[
A_\sigma^*A_{\sigma}(i,j)=\left\{
                    \begin{array}{cc}
                      \sum_{k\notin\sigma} \bar a_{k,i}a_{k,j}, & i,j\in\sigma \\
                      \sum_{k\in\sigma} \bar a_{k,i}a_{k,j}, & i,j\notin\sigma  \\
                      0 &\mbox{otherwise}.\\
                    \end{array}
                    \right.
                  \]
It follows that for all $i$ and $j$,
\begin{multline*}
{{2m}\choose{m}}\rm{Ave}_{\sigma}A_\sigma^*A_{\sigma}(i,j)\\=2\sum_{\{\sigma;\ i,j\in\sigma\}}\sum_{k\notin\sigma}\bar a_{k,i}a_{k,j}=2\sum_{k\not=i,j}\#\{\sigma;\ i,j\in\sigma,\ k\notin\sigma\} \bar a_{k,i}a_{k,j}.
\end{multline*}
If $i=j$ and $k\not=i$, then $\#\{\sigma;\ i,j\in\sigma,\ k\notin\sigma\}={{2m-2}\choose{m-1}}$. If
$i\not=j$ and $k\not=i,j$, then $\#\{\sigma;\ i,j\in\sigma,\ k\notin\sigma\}={{2m-3}\choose{m-2}}$.
This translates to
\begin{equation}\label{eq:ave}
\rm{Ave}_{\sigma}A_\sigma^*A_{\sigma}(i,j)=\left\{
  \begin{array}{lll}
   & 2{{2m-2}\choose{m-1}}{{2m}\choose{m}}^{-1}\sum_{k\not=i}|a_{k,i}|^2, & i=j \\
   & 2{{2m-3}\choose{m-2}}{{2m}\choose{m}}^{-1}\sum_{k\not=i,j}\bar a_{k,i}a_{k,j}, & i\not=j.
  \end{array}
\right.
\end{equation}
Since
\[
{{2m-2}\choose{m-1}}{{2m}\choose{m}}^{-1}=\frac{m}{4m-2}\ \ \ \ \mbox{and}\ \ \ \
{{2m-3}\choose{m-2}}{{2m}\choose{m}}^{-1}=\frac{m}{8m-4},
\]
we get from (\ref{eq:conv}), (\ref{eq:ave}) and the fact that $A$ has zero diagonal that
\[
\rm{Ave}_{\sigma}\rm{Tr}(A_{\sigma}^*A_{\sigma})^{p/2}\ge
\rm{Tr}\left(\left(\frac{m}{4m-2}\right)A^*A+\left(\frac{m}{4m-2}\right)B\right)^{p/2},
\]
where $B$ is the diagonal matrix whose $i$-th element is $\sum_k|a_{k,i}|^2$.
Proposition 1 in \cite{pe} now implies that
\[
\rm{Ave}_{\sigma}\rm{Tr}(A_{\sigma}^*A_{\sigma})^{p/2}\ge
\left(\frac{m}{4m-2}\right)^{p/2}\rm{Tr}(A^*A)^{p/2}\ge
\frac{1}{2^p}\rm{Tr}(A^*A)^{p/2}.\quad\quad\quad\endpf
\]

\bigskip

\noindent{\bf Proof of Theorem \ref{thm:paving}}
We shall use the well known (and easy to prove by interpolation) inequality
\begin{equation}\label{eq:uc}
\left(\frac12(\|x+y\|_p^p+\|x-y\|_p^p)\right)^{1/p}\le (\|x\|_p^{p/(p-1)}+\|y\|_p^{p/(p-1)})^{(p-1)/p}
\end{equation}
for all $2\le p<\infty$ and all matrices $x,y$ (see e.g.
\cite[Th. 5.1]{px}).
By Proposition \ref{pr:Asigma}, we can find a $\sigma\subset\{1,2,\dots,2m\}$ of cardinality $m$
such that, putting $P$ to be the natural projection onto the span of $\{e_i\}_{i\in\sigma}$, $Q=I-P$,
and $v=PAQ+QAP$,
\[
\|v\|_p^p\ge\frac{1}{2^p}\|A\|_p^p.
\]
Put also $u=PAP+QAQ$ and notice that $A=u+v$ and that $\|u-v\|_p=\|u+v\|_p=\|A\|_p$. Indeed,
$P-Q$ is a unitary transformation, so   $\|u-v\|_p=\|(P-Q)A(P-Q)\|_p=\|A\|_p$. Apply now
(\ref{eq:uc}) to $x=u+v, y=u-v$, to get
\[
\left(\frac12(\|2u\|_p^p+\|2v\|_p^p)\right)^{1/p}\le (\|u+v\|_p^{p/(p-1)}+\|u-v\|_p^{p/(p-1)})^{(p-1)/p}
\]
or
\[
(\|u\|_p^p+\|v\|_p^p)^{1/p}\le\|A\|_p.
\]
Thus,
\[
\|u\|_p\le \left(1-\frac{1}{2^p}\right)^{1/p}\|A\|_p.
\]
\endpf

\medskip

\noindent {\bf Acknowledgement:} I wish to thank the referee for detecting inaccuracies in the original version of this paper and for helping to improve the presentation, especially in section \ref{sec:complemented}.

%
%

\noindent G. Schechtman\\
Department of Mathematics\\
Weizmann Institute of Science\\
Rehovot, Israel\\
{\tt gideon@weizmann.ac.il}

\end{document}